\documentclass[11pt]{amsart}

\usepackage{amsmath}
\usepackage{amssymb}

\newtheorem{theorem}{Theorem}[section]

\newtheorem{lemma}[theorem]{Lemma}

\theoremstyle{definition}
\newtheorem{definition}[theorem]{Definition}

\newtheorem{conjecture}[theorem]{Conjecture}

\theoremstyle{remark}
\newtheorem{remark}[theorem]{Remark}

\newtheorem{notation}[theorem]{Notation}

\newcount\skewfactor
\def\mathunderaccent#1#2 {\let\theaccent#1\skewfactor#2
\mathpalette\putaccentunder}
\def\putaccentunder#1#2{\oalign{$#1#2$\crcr\hidewidth
\vbox to.2ex{\hbox{$#1\skew\skewfactor\theaccent{}$}\vss}\hidewidth}}

\def\smallbox#1{\leavevmode\thinspace\hbox{\vrule\vtop{\vbox
   {\hrule\kern1pt\hbox{\vphantom{\tt/}\thinspace{\tt#1}\thinspace}}
   \kern1pt\hrule}\vrule}\thinspace}

\newcommand{\bool}{{\bf B}}

\newcommand{\cf}{{\rm cf}}
\newcommand{\Depth}{{\rm Depth}}

\newcommand{\st}{{such that}}
\newcommand{\seq}{{sequence}}
\newcommand{\cont}{{continuous}}
\newcommand{\incr}{{increasing}}
\newcommand{\Then}{{\underline{Then}}}

\newcommand{\mat}{\mathcal}

\def\qedref#1{$\qed_{\reforiginal{#1}}$}

\title{depth of Boolean Algebras}
\author{Shimon Garti}
\address{Institute of Mathematics
 The Hebrew University of Jerusalem
 Jerusalem 91904, Israel}
\email{shimonygarty@hotmail.com}
\thanks{ }

\author{Saharon Shelah}
\address{Institute of Mathematics
 The Hebrew University of Jerusalem
 Jerusalem 91904, Israel
 and  Department of Mathematics
 Rutgers University
 New Brunswick, NJ 08854, USA}
\email{shelah@math.huji.ac.il}
\urladdr{http://www.math.rutgers.edu/\char`\~shelah}
\thanks{First typed: April 2007 \\
Research supported by the United States-Israel Binational
Science Foundation. Publication 911 of the second author}


\begin{document}
\let\labeloriginal\label
\let\reforiginal\ref
\def\ref#1{\reforiginal{#1}}
\def\label#1{\labeloriginal{#1}}

\begin{abstract}

Suppose $D$ is an ultrafilter on $\kappa$ and $\lambda^\kappa = \lambda$. We prove that if $\bool_i$ is a
Boolean algebra for every $i < \kappa$ and $\lambda$
bounds the ${\rm Depth}$ of every $\bool_i$, then the ${\rm Depth}$ of the ultraproduct mod $D$ is bounded by $\lambda^+$.
We also show that for singular cardinals with small cofinality, there is no gap at all.
This gives a full answer
to this problem in the constructible universe.

\end{abstract}

\maketitle

% start document here:
\newpage

\section{introduction}

Let ${\bf B}$ be a Boolean Algebra.
We define the depth of it as the supremum
on the cardinalities of well-ordered subsets in ${\bf B}$.
Now suppose that $\langle {\bf B_i} : i < \kappa \rangle$ is
a sequence of Boolean algebras, and $D$ is an ultrafilter on $\kappa$.
Define the ultra-product algebra ${\bf B}$ as $\prod\limits_{i<\kappa} \bool_i /D$.
The question (raised also for other cardinal invariants,
by Monk, in \cite{MR1077622}) is about the relationship between
${\rm Depth} (\bool)$ and $\prod_{i<\kappa}
{\rm Depth} (\bool_i) / D$. \newline

\par \noindent
Let us try to draw the picture: \newline
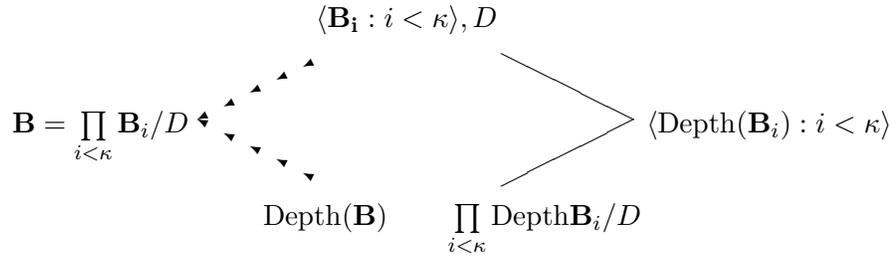
\begin{picture}(80,90)
\put(70,40){\line(4,2){50}}
\put(70,40){\line(4,-2){50}}
\put(235,40){\line(-2,1){50}}
\put(235,40){\line(-2,-1){50}}
\put(115,75){$\langle {\bf B_i} : i < \kappa \rangle, D$}
\put(95,0){${\rm Depth} (\bool)$}
\put(0,35){$\bool=\prod\limits_{i<\kappa} \bool_i /D$}
\put(240,35){$\langle {\rm Depth} (\bool_i) : i<\kappa \rangle$}
\put(165,0){$\prod\limits_{i<\kappa} {\rm Depth} \bool_i /D$}
\end{picture}

\medskip
\medskip

As we can see from the picture, given a sequence of Boolean algebras (of length $\kappa$) and an ultrafilter on $\kappa$, we have two alternating ways to produce a cardinal value. The left course creates, first, a new Boolean algebra namely the ultraproduct algebra $\bool$. Then we compute the Depth of it. In the second way, first of all we get rid of the algebraic structure, producing a sequence of cardinals (namely $\langle {\rm Depth} (\bool_i) : i<\kappa \rangle$). Then we compute the cardinality of its cartesian product divided by $D$. \newline

Shelah proved in \cite{MR1918108} \textsection 5, under the assumption $\bf {V} = \bf {L}$, that if $\lambda = \lambda^\kappa$ and
$\kappa = \cf (\kappa) < \lambda$, then you can build a sequence of
Boolean algebras $\langle {\bf B_i} : i < \kappa \rangle$,
such that ${\rm Depth} (\bool) > \prod_{i<\kappa}
{\rm Depth} (\bool_i) / D$ for every uniform ultrafilter $D$. This result is based on the square principle, introduced and proved in $\bf {L}$ by Jensen. \newline
A natural question is how far can this gap reach. We prove that
if $\bf {V} = \bf {L}$ then the gap is at most one cardinal.
In other words, for every regular cardinal and for every singular cardinal with high cofinality
we can create a gap (having the
square for every infinite cardinal in ${\bf L}$), but it
is limited to one cardinal. \newline
Observe that the assumption $\bf {V} = \bf {L}$ is just to make sure that every ultrafilter is regular.
We observe also that by \cite{MR2217966}, under some reasonable assumptions, there is no gap at all above a compact cardinal.\newline
We can ask further what happens if $\cf (\lambda) < \lambda$, and $\kappa \geq \cf (\lambda)$.
We prove here that if $\lambda$ is singular with small
cofinality, (i.e., all the cases which are not covered in the previous paragraph), then $\prod_{i<\kappa}
{\rm Depth} (\bool_i) / D \geq {\rm Depth} (\bool)$.
It is interesting to know that similar result holds
above a compact cardinal for singular cardinals with countable cofinaliy.
We suspect that it holds (for such cardinals) in ZFC. \newline
The proof of those results is based on an improvement to the
main Theorem in \cite{temp878}. It says that under some assumptions we can dominate the gap between ${\rm Depth} (\bool)$ and $\prod_{i<\kappa}
{\rm Depth} (\bool_i) / D$. In this paper we use weaker assumptions. We give here the full proof, so the
paper is self-contained. We intend to shed light on the other side of the coin (i.e., under large cardinals assumptions) in a subsequent paper.

\newpage

\section{The main theorem}

\begin{definition}
\label{0.1}
$\Depth$. \newline
Let ${\bf B}$ be a Boolean Algebra.
$$
\Depth(\bool):= \sup\{\theta:\exists \bar b=(b_\gamma:\gamma<\theta),
\hbox{ \incr\ \seq\ in } \bool\}
$$
\end{definition}

We use also an important variant of the Depth:

\begin{definition}
\label{0.2}
$\Depth^+$. \newline
Let $\bool$ be a Boolean Algebra.
$$
\Depth^+ (\bool):=\sup \{\theta^+:\exists \bar b=(b_\gamma:\gamma<\theta),
\hbox{ \incr\ \seq\ in } \bool\}
$$
\end{definition}

\par \noindent
Through the paper, we use the following notation:

\begin{notation}
\label{1.1}
\begin{enumerate}
\item[(a)] $\kappa,\lambda$ are infinite cardinals
\item[(b)] $D$ is a uniform ultrafilter on $\kappa$
\item[(c)] $\bool_i$ is a Boolean Algebra, for any $i<\kappa$
\item[(d)] $\bool=\prod\limits_{i<\kappa} \bool_i /D$
\item[(e)] for $\kappa < \lambda$, $S^\lambda_\kappa = \{\alpha < \lambda : \cf (\alpha) = \kappa \}$.
\end{enumerate}
\end{notation}

We state our main result: \newline

\begin{theorem}
\label{1.2}
Assume
\begin{enumerate}
\item[(a)] $\lambda=\cf(\lambda)$
\item[(b)] $\lambda = \lambda^\kappa$
\item[(c)] $\Depth^+ (\bool_i) \leq \lambda$, for every $i<\kappa$.
\end{enumerate}
\Then\ $\Depth^+ (\bool) \leq \lambda^+$.
\end{theorem}

\par \noindent
\emph {Proof}. \newline
Assume towards a contradiction that $\langle a_\alpha:\alpha<\lambda^+ \rangle$ is an \incr\ \seq\ in $\bool$. Let us write $a_\alpha$ as
$\langle a^\alpha_i:i<\kappa\rangle/D$ for every $\alpha<\lambda^+$.

Let $\langle M_\alpha:\alpha<\lambda^+\rangle$ be a continuous and increasing \seq\ of elementary submodels of $({\mat H} (\chi),\in)$ for sufficiently large
$\chi$ with the following properties $(\forall \alpha <\lambda^+)$: \newline
$(*)_1$
\begin{enumerate}
\item[(a)] $\| M_\alpha\|=\lambda$
\item[(b)] $\lambda+1 \subseteq M_\alpha$
\item[(c)] $\langle M_\beta:\beta \leq \alpha \rangle \in M_{\alpha+1}$
\item[(d)] $[M_\alpha]^\kappa \subseteq M_{\alpha+1}$.
\end{enumerate}

We may assume that $\langle a^\alpha_i:
\alpha<\lambda^+, i<\kappa\rangle \in M_0$.
We also assume that $\bool, \langle \bool_i : i<\kappa \rangle, D \in M_0$.

We will try to create a set $Z$, in the Lemma below, with
the following properties: \newline
$(*)_2$
\begin{enumerate}
\item[(a)] $Z \subseteq \lambda^+, |Z|=\lambda$
\item[(b)] $\exists i_* \in \kappa$ \st\ for every $\alpha<\beta,
\alpha,\beta\in Z$,
we have $\bool_{i_*} \models a^\alpha_{i_*}
<a^\beta_{i_*}$
\end{enumerate}
Since $|Z|=\lambda$, we have an \incr\ \seq\ of length $\lambda$
in $\bool_{i_*}$, so $\Depth^+ (\bool_{i_*}) \geq \lambda^+$, contradicting
the assumptions of the Theorem. \newline
\medskip
\hfill \qedref{1.2}

\begin{lemma}
\label{1.5}
There exists $Z$ as above.
\end{lemma}

\par \noindent
\emph {Proof}. \newline
For every $\alpha<\beta<\lambda^+$, define:
$$
A_{\alpha,\beta}=\{i<\kappa:\bool_i \models a^\alpha_i<a^\beta_i\}
$$
By the assumption, $A_{\alpha,\beta} \in D$
for all $\alpha<\beta<\lambda^+$.
Define $C:= \{ \gamma < \lambda^+ : \gamma = M_\gamma \cap \lambda^+ \}$, and $S := C \cap S^{\lambda^+}_\lambda$. Since $C$ is a club subset of $\lambda^+$, $S$ is a stationary subset of $\lambda^+$. Choose $\delta^*$ as the $\lambda$-th member of $S$.
For every $\alpha<\delta^*$,
Let $A_\alpha$ denote the set $A_{\alpha,\delta^*}$. \newline

Let $u\subseteq \delta^*$, $|u|\leq \kappa$. Notice that $u \in M_{\delta^*}$, by (d) of $(*)_1$ above. Define:
$$
S_u=\{\beta<\delta^*:\beta>{\rm sup} (u) \hbox{ and }
(\forall \alpha\in u) (A_{\alpha,\beta}=A_\alpha)\}.
$$

Choose $\delta_0 = 0$.
Choose $\delta_{\epsilon+1}\in S$ for every $\epsilon<\lambda$ \st\
$\epsilon<\zeta\Rightarrow {\rm sup}
\{\delta_{\epsilon+1}:\epsilon<\zeta\}<\delta_{\zeta+1}$.
Define $\delta_\epsilon$ to be the
limit of $\delta_{\gamma+1}$, when $\gamma<\epsilon$,
for every limit $\epsilon<\lambda$. Since ${\rm otp}(S \cap \delta^*) = \lambda$, we have:
\begin{enumerate}
\item[(a)] $\langle \delta_\epsilon:\epsilon<\lambda\rangle$
is \incr\ and \cont\
\item[(b)] $\rm sup \{\delta_\epsilon:\epsilon<\lambda\} = \delta^*$
\item[(c)] $\delta_{\epsilon+1} \in S$, for every $\epsilon<\lambda$
\end{enumerate}
Define, for every $\epsilon<\lambda$, the following family:
$$
{\frak A}_\epsilon=\{S_u\cap \delta_{\epsilon+1}\setminus
\delta_\epsilon:u\in [{\delta_{\epsilon+1}}]^{\leq \kappa}\}.
$$
We get a family of non-empty sets, which is downward $\kappa^+$-directed.
So, there is a $\kappa^+$-complete filter $E_\epsilon$ on
$[\delta_\epsilon, \delta_{\epsilon+1})$, with ${\frak A}_\epsilon
\subseteq E_\epsilon$, for every $\epsilon<\lambda$.

Define, for any $i<\kappa$ and $\epsilon<\lambda$, the sets
$W_{\epsilon,i} \subseteq [\delta_\epsilon,\delta_{\epsilon+1})$
and $B_\epsilon\subseteq \kappa$, by:
$$
W_{\epsilon,i} := \{\beta:\delta_\epsilon\leq \beta<\delta_{\epsilon+1}
\hbox{ and } i\in A_{\beta, \delta_{\epsilon+1}}\}
$$
$$
B_\epsilon:=\{i<\kappa:W_{\epsilon,i}\in E^+_\epsilon\}.
$$
Finally, take a look at $W_\epsilon:=\cap \{[\delta_\epsilon,
\delta_{\epsilon+1})\setminus W_{\epsilon,i}:
i\in \kappa \setminus B_\epsilon\}$.
For every $\epsilon<\lambda, W_\epsilon\in E_\epsilon$, since
$E_\epsilon$ is $\kappa^+$-complete, so clearly $W_\epsilon\neq \emptyset$.

Choose
$\beta=\beta_\epsilon\in W_\epsilon$.
If $i\in A_{\beta,\delta_{\epsilon+1}}$,
then $W_{\epsilon,i} \in E^+_\epsilon$, so $A_{\beta,\delta_{\epsilon+1}}\subseteq
B_\epsilon$ (by the definition of $B_\epsilon$). But,
$A_{\beta,\delta_{\epsilon+1}}\in D$, so
$B_\epsilon\in D$. For every $\epsilon < \lambda$, $A_{\delta_{\epsilon+1}}$ (which is $A_{\delta_{\epsilon+1}, \delta^*}$) belongs to $D$, so $B_\epsilon \cap A_{\delta_{\epsilon+1}} \in D$.
\newline
Choose $i_\epsilon\in B_\epsilon \cap A_{\delta_{\epsilon+1}}$, for every $\epsilon<\lambda$.
You have chose $\lambda\  i_\epsilon$-s from $\kappa$, so we
can arrange a fixed $i_*\in \kappa$ \st\ the set $Y=\{\epsilon<\lambda:\epsilon$
is an even ordinal, and $i_\epsilon=i_*\}$ has cardinality $\lambda$.

The last step will be as follows: \newline
define $Z=\{\delta_{\epsilon+1}:\epsilon\in Y\}$.
Clearly, $Z\in [\delta^*]^\lambda \subseteq [\lambda^+]^\lambda$.
We will show that for $\alpha<\beta$ from $Z$ we get
$\bool_{i_*} \models a^\alpha_{i_*}<a^\beta_{i_*}$.
The idea is that if $\alpha<\beta$ and $\alpha,\beta\in Z$,
then $i_*\in A_{\alpha,\beta}$.

Why? Recall that $\alpha=\delta_{\epsilon+1}$ and
$\beta=\delta_{\zeta+1}$, for some $\epsilon<\zeta<\lambda$
(that's the form of the members of $Z$). Define: \newline
$U_1 := S_{\{\delta_{\epsilon+1}\}} \cap
[\delta_\zeta,\delta_{\zeta+1}) \in
{\frak A}_\zeta \subseteq E_\zeta$. \newline
$U_2 := \{\gamma:\delta_\zeta\leq \gamma<\delta_{\zeta+1}$ and
$i_*\in A_{\gamma,\delta_{\zeta+1}}\}\in E^+_\zeta$. \newline
So, $U_1\cap U_2\neq \emptyset$. \newline
Choose $\iota\in U_1\cap U_2$.\newline
Now the following statements hold:
\begin{enumerate}
\item[(a)] $\bool_{i_*} \models a^\alpha_{i_*}<a^\iota_{i_*}$

[Why? Well, $\iota\in U_1$, so $A_{\delta_{\epsilon+1,\iota}}=
A_{\delta_{\epsilon+1}}$. But, $i_*\in B_\epsilon \cap A_{\delta_{\epsilon+1}} \subseteq A_{\delta_{\epsilon+1}}$, so
$i_*\in A_{\delta_{\epsilon+1,\iota}}$, which means that
$\bool_{i_*} \models a^{\delta_{\epsilon+1}}_{i_*}
 (=a^\alpha_{i_*})<a^\iota_{i_*}$].
\item[(b)] $\bool_{i_*}\models a^\iota_{i_*}<a^\beta_{i_*}$

[Why? Well, $\iota\in U_2$, so $i_*\in A_{\iota,\delta_{\zeta+1}}$,
which means that $\bool_{i_*} \models a^\iota_{i_*}<
a^{\delta_{\zeta+1}}_{i_*} (=a^\beta_{i_*})$].
\item[(c)] $\bool_{i_*}\models a^\alpha_{i_*}<a^\beta_{i_*}$

[Why? By (a)+(b)].
\end{enumerate}
So, we are done.

\hfill \qedref{1.5}

\newpage

\section{Depth in $\bf L$}

As a consequence of the main result from the previous section we have, under the constructibility axiom, as follows:

\begin{theorem}
\label{inl}
($GCH$) \newline
Assume
\begin{enumerate}
\item[(a)] $\kappa < \lambda$
\item[(b)] $\Depth (\bool_i) \leq \lambda$, for every $i<\kappa$.
\item[(c)] $\lambda = {\rm lim}_D (\langle \Depth (\bool_i) : i< \kappa \rangle )$
\end{enumerate}
\Then\ $\Depth (\bool) \leq \lambda^+$.
\end{theorem}

\par \noindent
\emph {Proof}. \newline
For every successor cardinal $\lambda^+$ we have (under the GCH)
$$
(\lambda^+)^\kappa=(2^\lambda)^\kappa=2^{\lambda \cdot \kappa}=2^\lambda=\lambda^+
$$
Clearly, $\lambda^+$ is a regular cardinal, and by (b) we know that
$\Depth^+ (\bool_i) \leq \lambda^+$ for every $i < \kappa$.
Now apply Theorem \ref{1.2} and conclude that $\Depth^+ (\bool) \leq \lambda^{+2}$, so $\Depth (\bool) \leq \lambda^+$ as required.

\hfill \qedref{inl}

\begin{remark}
\label{gch}
In $\bf L$ equality holds. The proof is similar to the proof in Theorem \ref{singu} below.
\end{remark}

\medskip

So if $\lambda$ is regular and $\kappa < \lambda$, or even $\lambda > \cf(\lambda) > \kappa$, we can build in $\bf L$ an example for $\Depth (\bool) > \prod_{i<\kappa} {\rm Depth} (\bool_i) / D$, but the discrepancy is just one cardinal. We can ask what happens if $\lambda$ is singular with small cofinality. The following Theorem gives an answer. Notice that this answers problem No. 12 from \cite{MR1393943}, for the case of singular cardinals with countable cofinality.

\begin{theorem}
\label{singu}
($V = L$) \newline
Assume
\begin{enumerate}
\item[(a)] $\cf (\lambda) < \lambda$
\item[(b)] $\kappa \geq \cf (\lambda)$
\item[(c)] $\Depth (\bool_i) \leq \lambda$, for every $i<\kappa$
\item[(d)] $\lambda = {\rm lim}_D (\langle \Depth (\bool_i) : i< \kappa \rangle )$.
\end{enumerate}
\Then\ $\Depth (\bool) = \prod_{i<\kappa}
{\rm Depth} (\bool_i) / D$.
\end{theorem}

\par \noindent
\emph {Proof}. \newline
First we claim that $\prod_{i<\kappa} {\rm Depth} (\bool_i) / D = \lambda^+$.
The basic idea is that in $\bf L$ we know that $D$ is regular (by the fundamental result of Donder, from \cite{MR969944}), so $\prod_{i<\kappa} {\rm Depth} (\bool_i) / D = \lambda^\kappa = \lambda^+$ (recall that $\cf (\lambda) \leq \kappa$). \newline
Now $\Depth (\bool) \geq \prod_{i<\kappa} {\rm Depth} (\bool_i) / D = \lambda^+$, by Theorem 4.14 from \cite{MR1077622} (since $\bf L$ $\models \rm GCH$).
On the other hand, Theorem \ref{inl} makes sure that $\Depth (\bool) \leq \lambda^+$ (by (c) of the present Theorem). So $\prod_{i<\kappa} {\rm Depth} (\bool_i) / D = \lambda^+ = \Depth (\bool)$, and we are done.

\hfill \qedref{singu}

\medskip

We know that if $\kappa$ is less than the first measurable cardinal, then every uniform ultrafilter on $\kappa$ is $\aleph_0$-regular. It gives us the result of
Theorem \ref{singu} for singular cardinals with countable cofinality, if the length of the sequence (i.e., $\kappa$) is below the first measurable. \newline
We have good evidence that something similar holds for singular cardinals with countable cofinality above a compact cardinal.
Moreover, if $\cf (\lambda) = \aleph_0$ then $\kappa \geq \cf (\lambda)$ for every infinite cardinal $\kappa$. It means that it is consistent with ZFC not to have a counterexample in this case.
So the following conjecture does make sense:

\begin{conjecture}
(ZFC) \newline
Assume
\begin{enumerate}
\item[(a)] $\aleph_0 = \cf (\lambda) < \lambda$, and $2^{\aleph_0} < \lambda$
\item[(b)] $\kappa < \lambda$
\item[(c)] $\Depth (\bool_i) \leq \lambda$, for every $i<\kappa$
\item[(d)] $\lambda = {\rm lim}_D (\langle \Depth (\bool_i) : i< \kappa \rangle )$.
\end{enumerate}
\Then\ $\Depth (\bool) \leq \prod_{i<\kappa}
{\rm Depth} (\bool_i) / D$.
\end{conjecture}

Notice that by \cite{MR1070264} we know that this question is independent when $2^{\aleph_0} > \lambda$, as
follows from Theorem 3.2 there.

\bibliographystyle{amsplain}
\bibliography{arlist}

\providecommand{\bysame}{\leavevmode\hbox to3em{\hrulefill}\thinspace}
\providecommand{\MR}{\relax\ifhmode\unskip\space\fi MR }
\providecommand{\MRhref}[2]{%
  \href{http://www.ams.org/mathscinet-getitem?mr=#1}{#2}
}
\providecommand{\href}[2]{#2}
\begin{thebibliography}{1}

\bibitem{MR969944}
Hans-Dieter Donder, \emph{Regularity of ultrafilters and the core model},
  Israel J. Math. \textbf{63} (1988), no.~3, 289--322. \MR{MR969944
  (90a:03071)}

\bibitem{temp878}
Shimon Garti and Saharon Shelah, \emph{On depth and {$depth\sp +$} of {B}oolean
  algebras}, Algebra Universalis, accepted.

\bibitem{MR1077622}
J.~Donald Monk, \emph{Cardinal functions on {B}oolean algebras}, Lectures in
  Mathematics ETH Z\"urich, Birkh\"auser Verlag, Basel, 1990. \MR{MR1077622
  (92d:06033)}

\bibitem{MR1393943}
\bysame, \emph{Cardinal invariants on {B}oolean algebras}, Progress in
  Mathematics, vol. 142, Birkh\"auser Verlag, Basel, 1996. \MR{MR1393943
  (97c:06018)}

\bibitem{MR1070264}
Saharon Shelah, \emph{Products of regular cardinals and cardinal invariants of
  products of {B}oolean algebras}, Israel J. Math. \textbf{70} (1990), no.~2,
  129--187. \MR{MR1070264 (91i:03102)}

\bibitem{MR1918108}
\bysame, \emph{More constructions for {B}oolean algebras}, Arch. Math. Logic
  \textbf{41} (2002), no.~5, 401--441. \MR{MR1918108 (2003f:03063)}

\bibitem{MR2217966}
\bysame, \emph{The depth of ultraproducts of {B}oolean algebras}, Algebra
  Universalis \textbf{54} (2005), no.~1, 91--96. \MR{MR2217966 (2007b:06017)}

\end{thebibliography}

\end{document}